\documentclass[11pt]{article}
\usepackage{lipsum}

\usepackage{epsfig}
\usepackage{amssymb}
\usepackage{color}
\usepackage{enumerate}
\usepackage{xmpmulti}
\usepackage{amssymb}
\usepackage{rotating}
\usepackage{amsmath}
\usepackage{amsthm}
\usepackage{amscd}
\usepackage[square,comma,sort&compress,numbers]{natbib} 
\usepackage{epsfig}
\usepackage{color}
\usepackage[latin5]{inputenc}
\usepackage{psfrag}
\usepackage{graphicx}

\newcommand{\cL}{{{\mathcal{L}}}}
\newcommand{\cC}{{{\mathcal{C}}}}

\newcommand{\MCG}{{\operatorname{MCG}}}

\newtheorem{thm}{Theorem}
\newtheorem{lem}[thm]{Lemma}

\theoremstyle{definition}		
\newtheorem{defn}[thm]{Definition}
\newtheorem{example}[thm]{Example}

\theoremstyle{remark}
\newtheorem{remark}[thm]{Remark}

\newtheorem*{prooftheorem4}{Proof of Theorem 4}

\newcommand\R{\mathbb{R}}

\newcommand\cA{{\mathcal A}}
\newcommand\cF{{\mathcal F}}
\newcommand\cS{{\mathcal S}}

\newcommand{\tropical}[1]{\left[\,#1\,\right]}

\newcommand{\I}{^{-1}}
\DeclareMathOperator{\PIL}{PIL}
\DeclareMathOperator{\PIP}{PIP}

\DeclareMathOperator{\MF}{\mathcal{MF}}

\newtheorem*{problem2}{\textnormal{\bf Problem 2}}
\newtheorem*{problem1}{\textnormal{\bf Problem 1}}

\graphicspath{{./}{figures/}}

\title{Applications of the Dynnikov coordinate system  on the boundary of Teichm\"uller space }
\author{S. \"Oyk\"u Yurtta\c s\footnote{Dicle University, Science
    Faculty, Mathematics Department, 21280, Diyarbak\i r, Turkey, e-mail: oykuyurttas@gmail.com}}

\begin{document}
\maketitle

\begin{abstract}

 The Dynnikov coordinate system puts global coordinates on the boundary of Teichm\"uller space of an $n$--punctured disk.  We survey the  Dynnikov coordinate system, and  investigate  how we use this coordinate system to study pseudo--Anosov braids  making use of results from Thurston's theory on surface homeomorphisms.  \end{abstract}
  \section{Introduction}

 Let $D_n$ be an $n$-punctured disk ($n\geq 3$).   A simple closed curve on $D_n$ is \emph{inessential} if it bounds an unpunctured disk, once punctured disk or an unpunctured annulus. It is \emph{essential} otherwise.   An \emph{integral lamination}~$\cL$~on~$D_n$~is a non-empty union of mutually disjoint unoriented essential simple closed curves in $D_n$ up to isotopy. The usual way to coordinatize integral laminations  is to use either Dehn-Thurston coordinates or  train track coordinates \cite{penner, W88}.   An alternative way which is more suitable than such combinatorial descriptions, particularly for problems on the finitely punctured disk, is to use Dynnikov's coordinates. 
 
 	The Dynnikov coordinate system gives a homeomorphism from the space of measured foliations $\mathcal{MF}_n$ (up to isotopy and Whitehead equivalence) on $D_n$ to $\mathbb{R}^{2n-4}\setminus \{0\}$; and restricts to a bijection from the set of integral laminations  $\mathcal{L}_n$ on $D_n$ to $\mathbb{Z}^{2n-4}\setminus \{0\}$  \cite{D02, wiest, M06, D08, or08, paper1, paper2, paper3, paper4}.  The mapping class group $\MCG(D_n)$ of $D_n$ is isomorphic to the $n$--braid group $B_n$ modulo its center~\cite{emil2} so that isotopy classes of orientation preserving homeomorphisms on $D_n$ can be represented by $n$--braids.  The action of $B_n$ on the space of Dynnikov coordinates is given by the \emph{update rules} of Theorem \ref{thm:update} below \cite{D02, wiest, M06, D08, or08, paper1, paper2, paper3,paper4}. 
	
	Dynnikov coordinate system together with the update rules have been used to solve many interesting problems such as giving a method for a solution of the word problem of $B_n$ \cite{or08, D08}, computing the topological entropy of braids \cite{M06, thiffeault1, paper1},  studying the dynamics of pseudo--Anosov braids \cite{paper3}, and  have a wide range of dynamical applications \cite{thiffeault1, thiffeault2, thiffeault3}. Furthermore, a recent application of the  Dynnikov coordinate system  \cite{paper4}  solves for surfaces of genus zero, a long-standing conjecture which asks the existence of a polynomial time algorithm that decides whether an integral lamination, specified in terms of a coordinate system, is connected or not. Namely,  in \cite{paper4} a quadratic time algorithm for calculating the number of components of an integral lamination from its Dynnikov coordinates is introduced.  A related problem is solved in \cite{paper6} where an algorithm  for calculating the geometric intersection number of two integral laminations on $D_n$ taking as input their Dynnikov coordinates is described. The algorithm has complexity polynomial in the number of punctures of $D_n$, and the sum of the absolute values of the Dynnikov coordinates of the integral laminations.
 
	In this survey we shall first investigate the Dynnikov coordinate system and the update rules defined for Artin braid generators. Then we shall  present an efficient method for studying pseudo-Anosov braids on $D_n$ which is based on results from Thurston's theory on surface homeomorphisms, and makes use of Dynnikov's coordinates together with the update rules on the boundary of Teichm\"uller space.

\section{The Dynnikov coordinate system}\label{sec:dynnikov}
\begin{figure}[htbp]
\begin{center}
\psfrag{1}[tl]{\begin{turn}{-90}$\scriptstyle{\alpha_1}$\end{turn}} 
\psfrag{b1}[tl]{\begin{turn}{-90}$\scriptstyle{\beta_1}$\end{turn}} 
\psfrag{b2}[tl]{\begin{turn}{-90}$\scriptstyle{\beta_{i+1}}$\end{turn}} 
\psfrag{bi}[tl]{\begin{turn}{-90}$\scriptstyle{\beta_{i}}$\end{turn}}
\psfrag{b6}[tl]{\begin{turn}{-90}$\scriptstyle{\beta_{n-1}}$\end{turn}} 
\psfrag{2i}[tl]{\begin{turn}{-90}$\scriptstyle{\alpha_{2i+2}}$\end{turn}}
\psfrag{2}[tl]{\begin{turn}{-90}$\scriptstyle{\alpha_2}$\end{turn}} 
\psfrag{2i5}[tl]{\begin{turn}{-90}$\scriptstyle{\alpha_{2i-3}}$\end{turn}}
\psfrag{2i2}[tl]{\begin{turn}{-90}$\scriptstyle{\alpha_{2i}}$\end{turn}}
\psfrag{2i3}[tl]{\begin{turn}{-90}$\scriptstyle{\alpha_{2i-1}}$\end{turn}}
\psfrag{2i1}[tl]{\begin{turn}{-90}$\scriptstyle{\alpha_{2i+1}}$\end{turn}}
\psfrag{2i4}[tl]{\begin{turn}{-90}$\scriptstyle{\alpha_{2i-2}}$\end{turn}}
\psfrag{n4}[tl]{\begin{turn}{-90}$\scriptstyle{\alpha_{2n-4}}$\end{turn}} 
\psfrag{n5}[tl]{\begin{turn}{-90}$\scriptstyle{\alpha_{2n-5}}$\end{turn}}
\psfrag{a}{$\scriptstyle{1}$} 
\psfrag{b}{$\scriptstyle{2}$} 
\psfrag{c}{$\scriptstyle{i}$} 
\psfrag{d}{$\scriptstyle{i+1}$} 
\psfrag{e}[r]{$\scriptstyle{i+2}$} 
\psfrag{f}{$\scriptstyle{n-1}$} 
\psfrag{g}{$\scriptstyle{n}$} 
\includegraphics[width=0.74\textwidth]{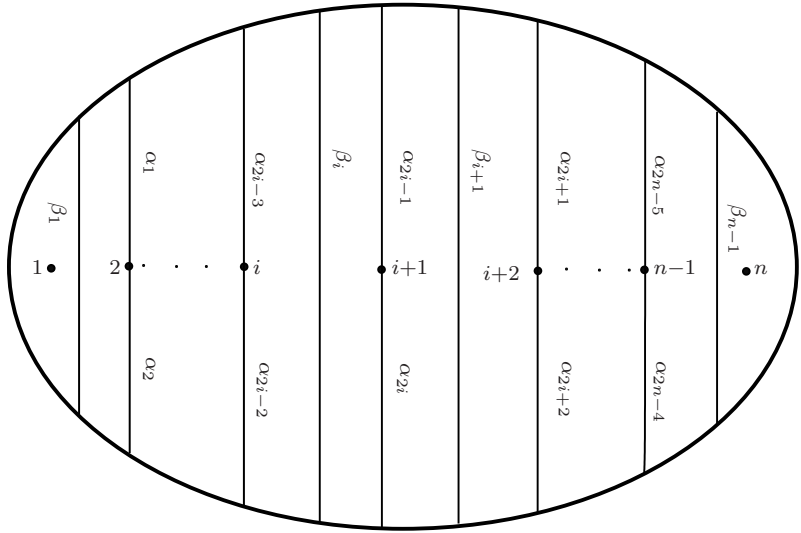}
\caption{The arcs $\alpha_i$ and $\beta_i$}\label{fig:dynn-arcs}
\end{center}
\end{figure}

Consider the Dynnikov arcs $\alpha_i $ ($1\le i\le 2n-4$) and $\beta_i$ ($1\leq i\leq n-1$) in~$D_n$~which have each endpoint either on the boundary of $D_n$ or at a puncture as shown in Figure~\ref{fig:dynn-arcs}. Each complementary region in Figure~\ref{fig:dynn-arcs} is triangular other than the two end
 regions (each region on the left and right side of the $i^{\text{th}}$ puncture for $2\leq i \leq n-1$ is a triangle since it is bounded by three arcs when the boundary of the disk is identified with a point), and there are $2n-4$ such triangles.  

\noindent Let $L$ be a minimal representative of an integral lamination $\mathcal{L}\in \mathcal{L}_n$ (i.e. $L$ intersects each Dynnikov arc minimally).  Let $\alpha_i$, $\beta_i$ denote the number of intersections of $L$ with the arc $\alpha_i$ and $\beta_i$ respectively.  It will always be clear from the context whether these symbols refer to arcs or intersection numbers assigned on the arcs. We call $(\alpha; \beta)=(\alpha_1,\dots,\alpha_{2n-4} ;\, \beta_1,\dots, \beta_{n-1})$  the {\em triangle coordinates} of $\mathcal{L}$.  For each $i$ with $1\leq i\leq n-2$, let  $S_i$ denote the subset of~$D_n$ bounded by the arcs $\beta_i$ and $\beta_{i+1}$.  

\begin{figure}[h!]
\centering
\psfrag{2i-1}[tl]{$\scriptstyle{\alpha_{2i-1}}$} 
\psfrag{i+1}[tl]{$\scriptstyle{\beta_{i+1}}$} 
\psfrag{i}[tl]{$\scriptstyle{\beta_{i}}$} 
\psfrag{2i}[tl]{$\scriptstyle{\alpha_{2i}}$}
\psfrag{d1}[tl]{$\scriptstyle{\Delta_{2i-1}}$} 
\psfrag{d2}[tl]{$\scriptstyle{\Delta_{2i}}$} 
\psfrag{si}[tl]{$\scriptstyle{S_{i}}$} 
\includegraphics[width=0.48\textwidth]{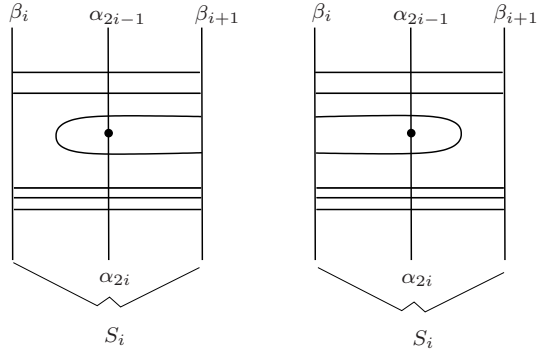}
\caption{Above, below, left loop and right loop components in $S_i$}\label{pathcomponentssurvey}
\end{figure}

A \emph{path component} of $L$ in $S_i$ is a component of $L\cap S_i$. By the minimal intersection condition, there are four types of path components of $L$ in each $S_i$:  an \emph{above component} has end points on $\beta_i$ and $\beta_{i+1}$ and intersects $\alpha_{2i-1}$ but not $\alpha_{2i}$. A \emph{below component} has end points on $\beta_i$ and $\beta_{i+1}$ and intersects $\alpha_{2i}$ but not  $\alpha_{2i-1}$.  A \emph{left loop component} has both end points on $\beta_{i+1}$ and intersects both of  $\alpha_{2i-1}$ and $\alpha_{2i}$.  A \emph{right loop component} has both end points on $\beta_{i}$, and intersects both of the arcs $\alpha_{2i-1}$ and $\alpha_{2i}$.  Figure \ref{pathcomponentssurvey} illustrates such path components. Clearly, there could only be one of the two types of loop components in region $S_i$  for each $1\leq i \leq n-2$ since the curves in $L$ are mutually disjoint. Define  
 \begin{align}\label{bi}
 b_i=\frac{\beta_i-\beta_{i+1}}{2}.
 \end{align}

  We immediately observe from  (\ref{bi}) that there are~$|b_i|$ loop components in each $S_i$, and when $b_i>0$ the loop components are right and when $b_i<0$ the loop components are left.  Note that there is one type of path component in the end regions: left loop components in the leftmost  region, and right loop components in the rightmost region. It follows that the number of loop components in the leftmost region equals $\frac{\beta_1}{2}$, and the number of loop components in the rightmost region equals $\frac{\beta_{n-1}}{2}$.  Furthermore, the numbers of above and below components in each $S_i$ are given by $\alpha_{2i-1}-|b_i|$ and $\alpha_{2i}-|b_i|$ respectively.  
  
    In Example \ref{const2} we shall reconstruct $\cL$ from its triangle coordinates. Namely, given the triangle coordinates of an integral lamination $\cL$, we shall determine the number of loop, above and below components in each region $S_i$, glue together these path components in a unique way up to isotopy,  and hence  construct ~$\cL$ uniquely.  
  
\begin{figure}[h!]
\centering
\psfrag{4}[tl]{$\scriptscriptstyle{4}$} 
\psfrag{2}[tl]{$\scriptscriptstyle{2}$} 
\psfrag{3}[tl]{$\scriptscriptstyle{3}$} 
\psfrag{1}[tl]{$\scriptscriptstyle{1}$} 
\psfrag{0}[tl]{$\scriptscriptstyle{0}$} 
\includegraphics[width=0.75\textwidth]{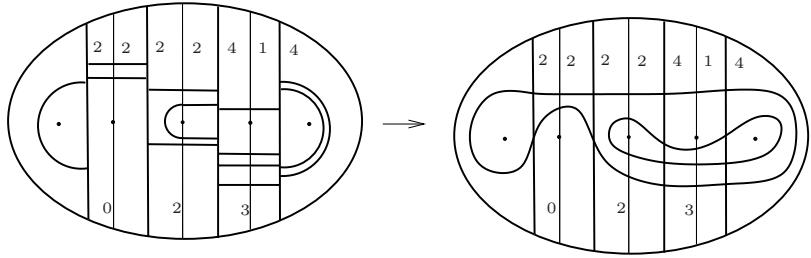}
\caption{To construct $\mathcal{L}$ from  triangle coordinates, we glue together above, below and loop components in each $S_i$ uniquely up to isotopy.}\label{const1}
\end{figure}

  \begin{example}\label{const2}
  Let $(2,0,2,2,1,3;\, 2,2,4,4)$ be the triangle coordinates of an integral lamination $\mathcal{L}\in \mathcal{L}_5$. First, we compute the number of loop components, and then the number of above and below components  in each region $S_i$. By (\ref{bi}) we compute that $b_1=0,~b_2=-1,~b_3=0$, and hence there are no loop components in $S_1$ and $S_3$, and  one left loop component in $S_2$. Also, there is $\frac{\beta_1}{2}=1$ loop component in the leftmost region, and $\frac{\beta_{n-1}}{2}=2$ loop components in the rightmost region.   Next, we work out the number of above and below components in each $S_i$. Since $\alpha_1-|b_1|=2$ and $\alpha_2-|b_1|=0$, there is no below component and 2 above components in $S_1$. Similar computations give that there is  $1$ below and $1$ above component in $S_2$, and $3$ below and $1$ above  components in $S_3$.  The path components in each $S_i$ are connected in a unique way up to isotopy, and hence we get $\cL$ as depicted in  Figure~\ref{const1}.  \end{example}

\vspace{-1mm}

However, it is not always possible to construct an integral lamination from given triangle coordinates since they must satisfy  the triangle inequality in each of the
triangular region of Figure~\ref{fig:dynn-arcs} as well as additional conditions to ensure that the curve system has no component which is boundary or puncture parallel.   For example,  an attempt to construct an integral lamination $\mathcal{L}\in \mathcal{L}_3$ with $(\alpha_1, \alpha_2, \beta_1, \beta_2)=(1,1,2,4)$ 
would fail since the triangle inequality is not satisfied in the triangular region bounded by the arcs  $\alpha_1,\alpha_2, \beta_2$ ($\beta_2>\alpha_1+\alpha_2$). Next, we define the Dynnikov coordinate system which coordinatizes $\mathcal{L}_n$ bijectively and with the least number of coordinates.
  
\begin{defn}
The {\em Dynnikov coordinate function}
$\rho\colon\cL_n\to\mathbb{Z}^{2n-4}\setminus \{0\}$ is defined by 
\[
\rho(\cL) = (a;\,b) = (a_1,\ldots,a_{n-2};\,b_1,\ldots,b_{n-2}),
\]
where 
\begin{equation}
\label{eq:dynn-coords}
a_i = \frac{\alpha_{2i}-\alpha_{2i-1}}{2} \qquad\text{and}\qquad
b_i = \frac{\beta_i - \beta_{i+1}}{2} \qquad
\end{equation}
for $1\le i\le n-2$.

\end{defn}

For example, the Dynnikov coordinates of the integral lamination depicted in Figure \ref{const1}  are given by $\rho(\cL) = (-1, 0,1;\, 0,-1,0)$. We already know the geometric interpretation of $b_i$ coordinates, also observe that $2a_i$ gives the difference between the number of below components and the number of above components in  $S_i$.  

	Theorem \ref{lem:dynninvert} \cite{paper1, paper2, paper3, paper4, paper6} gives formulae that recovers the triangle coordinates (and hence the integral lamination~$\cL$) from the Dynnikov coordinates.

\begin{thm}[Inversion of Dynnikov coordinates]
\label{lem:dynninvert}
Let $(a;\,b) \in \mathbb{Z}^{2n-4}\setminus\left\{0\right\}$. Then $(a;\,b)$ is the Dynnikov coordinate of exactly one element $\cL$ of $\cL_n$, which has 
\begin{align}
\beta_i&=2\max_{1\leq k\leq n-2}\left[|a_k|+\max(b_k,0)+\sum_{j=1}^{k-1}b_j\right]-2\sum_{j=1}^{i-1}b_j 
\label{opopp}
\end{align}
\begin{align}
\alpha_i&= \left\{ \begin{array}{ll}
         (-1)^i a_{\lceil i/2 \rceil}+\frac{\beta_{\lceil i/2\rceil}}{2}& \mbox{if $b_{\lceil i/2\rceil} \geq 0$};\\
        (-1)^i a_{\lceil i/2\rceil}+\frac{\beta_{1+\lceil i/2\rceil}}{2}& \mbox{if $b_{\lceil i/2\rceil} \leq 0$}   
         \end{array} \right.
\label{op}
\end{align}

\noindent where $\lceil x\rceil$ denotes the smallest integer which is not less than~$x$.

\end{thm}
A detailed proof of this theorem can be found \cite{paper2}. Here, we briefly explain the key idea behind the proof of the theorem.  We first observe that all of the $\beta_i$ can be computed from the coordinates~$b_j$ if~$\beta_1$ is known  since $\beta_i = \beta_1 - 2\displaystyle\sum^{i-1}_{j=1} b_j$ by (\ref{bi}).  Let $m_i$ denote the smaller of the number of above and below components in $S_i$. From Figure \ref{pathcomponentssurvey} we observe that 
$ \beta_i=2\left[\left|a_i \right|+\max(b_i,0)+m_i\right].$ And hence $$ \beta_1=2\left[\left|a_i \right|+\max(b_i,0)+m_i+\displaystyle\sum^{i-1}_{j=1}b_j \right]\quad  \text{for}\quad 1\leq i \leq n-2.$$

\noindent The crucial observation is for some $1\leq i\leq n-1$ we must have $m_i=0$ since the integral lamination doesn't contain any boundary parallel curves. Therefore for some $1\leq i\leq n-1$ we have

$$
\beta_1=2\left[\left|a_i \right|+\max(b_i,0)+\displaystyle\sum^{i-1}_{j=1} b_j \right]
$$

\noindent and hence,

\begin{eqnarray*}
\beta_1=\max_{1\leq k \leq{n-2}}2\left[{\left|a_k \right|+\max(b_k,0)+\displaystyle\sum^{k-1}_{j=1}}b_j\right].
\end{eqnarray*}

$\alpha_i$ can then be easily deduced using the coordinates~$a_j$.

\section{Update Rules}

In this section we explain how to compute the action of the $n$--braid group $B_n$ on the set of Dynnikov coordinates $\cC_n=\mathbb{Z}^{2n-4}\setminus\left\{0\right\}$. To do this, we shall use the update rules  \cite{D02, wiest, M06, D08, or08, paper1, paper2, paper3, paper4} of Theorem \ref{thm:update} which describe the action of each Artin's braid generator $\sigma_i$, $\sigma^{-1}_i,$ $(1 \leq i\leq n-1)$ on $\mathcal{L}_n$ in terms of Dynnikov coordinates. That is, the update rules tell us~$\rho(\sigma_i(\cL))$ and $\rho(\sigma_i^{-1}(\cL))$ in terms of $\rho(\cL)$. Here~$\sigma_i$ denotes the counterclockwise interchange of the $i^{\text{th}}$ and $i+1{^\text{th}}$ punctures, and the notation $x^+$ denotes $\max(x,0)$. For computational convenience, we work in the {\em max-plus semiring} $(\R,\oplus,\otimes)$ where $a\oplus b = \max(a;\,b)$ and $a\otimes b = a+b$ (so the multiplicative identity is $0$). For simplicity, we use normal additive and multiplicative notation in the formulae and express them in square brackets to indicate that the operations should be interpreted in their max-plus sense. That is,
$[a+b]=\max(a;\,b)$, $[ab]=a+b$, $[a/b] = a-b$, and $[1]=0$, the multiplicative identity. 

We note that the difference between the update rules given here and those that appeared in~\cite{D02,M06} will be that here $B_n$ acts on $D_n$ whereas in the cited papers $B_n$ acts on the central $n$ punctures in $D_{n+2}$. Thus we give special formulae for the action of $\sigma_1, \sigma_{n-1}$ and their inverses. 

 Let $(a' ; b')$  and  $(a'' ; b'')$ denote the Dynnikov coordinates of the integral laminations $\sigma_i(\cL)$ and  $\sigma^{-1}_i(\cL)$ respectively. 
 \begin{thm}\label{thm:update}
Let~$(a;\,b)\in\cC_n$ and \mbox{$1\le i\le n-1$}. Then $a_j'=a_j $, $b_j'=b_j$, $a_j''=a_j $ and $b_j''=b_j$ except when
$j=i-1$ or $j=i$, and:
\begin{itemize}
\item if $i=1$ then
\begin{align*}
a_1' &= \tropical{ \frac{a_1b_1}{a_1+1+b_1} }, &
b_1' &= \tropical{ \frac{1+b_1}{a_1} }\\*
a_1'' &= \tropical{ \frac{1+a_1(1+b_1)}{b_1}  }, &
b_1'' &= \tropical{ a_1(1+b_1) }; 
\end{align*}

\item if $2\le i \le n-2$ then
\begin{align*}
a_{i-1}' &= \tropical{ a_{i-1}(1+b_{i-1})+a_ib_{i-1}  }, &
b_{i-1}' &= \tropical{ \frac{a_ib_{i-1}b_i}{a_{i-1}(1+b_{i-1})(1+b_i)
+ a_ib_{i-1}}  }\\*
a_i' &= \tropical{ \frac{a_{i-1}a_ib_i}{a_{i-1}(1+b_i)+a_i}  }, &
b_i' &= \tropical{ \frac{a_{i-1}(1+b_{i-1})(1+b_i) + a_ib_{i-1}}{a_i} };\\*
a_{i-1}'' &= \tropical{
  \frac{a_{i-1}a_i}{a_{i-1}b_{i-1}+a_i(1+b_{i-1})}  },  &
b_{i-1}'' &= 
\tropical{ \frac{a_{i-1}b_{i-1}b_i}{a_{i-1}b_{i-1}+a_i(1+b_{i-1})(1+b_i)} }, \\*
a_i'' &= \tropical{ \frac{a_{i-1}+a_i(1+b_i)}{b_i}  }, &
b_i'' &= \tropical{ \frac{a_{i-1}b_{i-1}+a_i(1+b_{i-1})(1+b_i)}{a_{i-1}} };  
\end{align*}
\item if $i=n-1$ then
\begin{align*}
a_{n-2}' &= \tropical{ a_{n-2}(1+b_{n-2})+b_{n-2} }, &
b_{n-2}' &= \tropical{ \frac{b_{n-2}}{a_{n-2}(1+b_{n-2})} }\\*
a_{n-2}'' &= \tropical{\frac{a_{n-2}}{a_{n-2}b_{n-2}+1+b_{n-2}} }, &
b_{n-2}'' &= \tropical{ \frac{a_{n-2}b_{n-2}}{1+b_{n-2}} }.
\end{align*}
\end{itemize}

\end{thm}

\begin{figure}[h!]
\centering
\psfrag{a}[tl]{$\scriptstyle{a}$} 
\psfrag{b}[tl]{$\scriptstyle{b}$} 
\psfrag{c}[tl]{$\scriptstyle{c}$} 
\psfrag{d}[tl]{$\scriptstyle{d}$}
\psfrag{e}[tl]{$\scriptstyle{e}$} 
\psfrag{f}[tl]{$\scriptstyle{f}$} 
\psfrag{si}[tl]{$\scriptstyle{S_{i}}$} 
\includegraphics[width=0.45\textwidth]{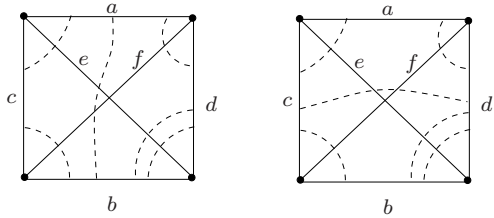}
\caption{Two possible cases for path components of $L$ in a quadrilateral}\label{tricksurvey}
\centering
\end{figure}

\noindent Before we prove Theorem \ref{thm:update} we give the following well known lemma \cite{dylan, D08} the proof of which is immediate from Figure \ref{tricksurvey}.

\begin{lem}\label{quadrilateral} 
Let $L$ be a minimal representative of an integral lamination $\mathcal{L}$, and $Q$ be a quadrilateral in $D_n$ with all of its vertices at punctures (where the boundary of $D_n$ is regarded as a puncture at $\infty$) and containing no punctures in its interior. Let $a, b, c, d,e,f$ denote the number of intersections of $L$ with the corresponding edges and diagonals of $Q$ as shown in Figure~\ref{tricksurvey}. Then,
\begin{eqnarray*}
e+f=\max(a+b, c+d)
\end{eqnarray*}
\end{lem}

\begin{figure}[h!]
\centering
\psfrag{u}[tl]{$\scriptstyle{u_{i}}$}
\psfrag{m3}[tl]{$\scriptstyle{\alpha _{2i-1}}$}
\psfrag{m2}[tl]{$\scriptstyle{\alpha _{2i-2}}$}
\psfrag{m4}[tl]{$\scriptstyle{\alpha_{2i}}$}
\psfrag{m3p}[tl]{$\scriptstyle{\alpha' _{2i-1}}$}
\psfrag{n3}[tl]{$\scriptstyle{\beta _{i+1}}$}
\includegraphics[width=0.55\textwidth]{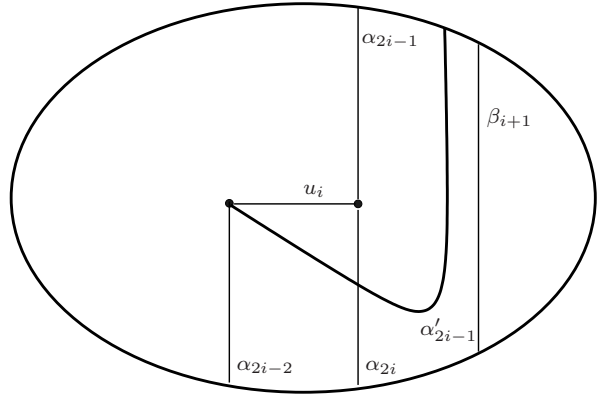}
\caption{A quadrilateral in $D_n$ with diagonal $\alpha'_{2i-1}$}\label{localaction3asurvey}
\end{figure}

\begin{prooftheorem4}
 Since $B_n$ acts on both the set of Dynnikov arcs  and $\mathcal{L}_n$, and the minimum intersection function is equivariant under this action,  the number of intersections of~$\sigma_i(L)$~with~the Dynnikov arcs $\alpha_j$  and $\beta_j$   equals the number of intersections of $\alpha'_j=\sigma^{-1}_i(\alpha_j)$ and $\beta'_j=\sigma^{-1}_i(\beta_j)$ with $L$.  Therefore we have,
$$a'_j=\frac{\alpha'_{2j}-\alpha'_{2j-1}}{2}  \quad\text{and}  \quad b'_j=\frac{\beta'_{j}-\beta'_{j+1}}{2}.$$
Set $A_{i}=2a_i$ and $B_{i}=2b_i$. Observe first that $\beta'_j=\beta_j$ for $j\neq i$ and $\alpha'_j=\alpha_j$ for $j<2i-3$ or $j>2i$. Therefore, $A'_j=A_j$ and $B'_j=B_j$ for except $j=i-1$ and $j=i$.   Here we compute $A'_i$.  The other coordinates are computed similarly.   Consider the quadrilateral in Figure \ref{localaction3asurvey} where $\alpha'_{2i-1}$ and $\alpha_{2i}$ are diagonals.  Then, $A'_i=\left[\frac{\alpha'_{2i}}{\alpha'_{2i-1}}\right]$. Let $u_i$ denote the number of intersections of $L$ with  the arc $u_i$ which has its end points on the $i^\text{th}$ and $i+1^\text{th}$ punctures as shown in Figure  \ref{localaction3asurvey}. By Lemma~\ref{tricksurvey} we get $$u_i+\beta_i=\max(\alpha_{2i-3}+\alpha_{2i},\alpha_{2i-2}+\alpha_{2i-1}).$$ That is,  $$u_i=\left[\frac{\alpha_{2i-3}\alpha_{2i}+\alpha_{2i-2}\alpha_{2i-1}}{\beta_i} \right].$$

\noindent  We have $\alpha'_{2i}=\alpha_{2i-2}$, and $\alpha'_{2i-1}+\alpha_{2i}=\max (u_i+\beta_{i+1},\alpha_{2i-2}+\alpha_{2i-1}).$  That is,

$$\alpha'_{2i-1}=\left[\frac{u_i\beta_{i+1}+\alpha_{2i-2}\alpha_{2i-1}}{\alpha_{2i}}\right].$$

Then we get, 
\begin{eqnarray*}
A'_i=\left[\frac{A_{i-1}A_iB_i}{A_{i-1}(1+B_i)+A_i}\right]
\end{eqnarray*}
That is,  $A'_i= A_{i-1} + A_i + B_i -\max(A_i, A_{i-1} + \max(0, B_i))$. Dividing both sides of the equation by $2$ replaces each  $A_i$ and $B_i$ with $a_i$ and $b_i$.

\end{prooftheorem4}

\section{Tools from Thurston's theory on surface homeomorphisms}
In this section we give  some basic terminology and principal results from Thurston's theory on surface homeomorphisms  \cite{FLP79, W88, BH95, mosher1, penner} which are necessary for our purposes. 

 Let $M$ be a surface of genus $g$ with  $b$ boundary components and $s$ punctures, which has negative Euler characteristic.  A \emph{measured foliation} $(\mathcal{F},\mu)$ on $M$ is a singular foliation $\mathcal{F}$ equipped with a measure $\mu$ such that  each arc $\alpha$ in $M$ transverse to $\mathcal{F}$ is assigned a positive number $\mu(\alpha)\in \R^{+} $ such that the following hold: 
\begin{itemize}
\item If $\alpha_1$ and $\alpha_2$ are two transverse arcs that are isotopic through other transverse arcs whose end points stay in the same leaves, then $\mu(\alpha_1)=\mu(\alpha_2).$
\item If $\alpha$ is a transverse arc such that $\alpha=\alpha_1\cup\alpha_2$ with $\alpha_1\cap\alpha_2=\partial\alpha_1\cap \partial\alpha_2$, then $\mu(\alpha)=\mu(\alpha_1)+\mu(\alpha_2).$
\end{itemize}

\begin{figure}[h!]
\centering
\psfrag{u}[tl]{$\scriptstyle{u_{i}}$}
\psfrag{m3}[tl]{$\scriptstyle{\alpha _{2i-1}}$}
\psfrag{m2}[tl]{$\scriptstyle{\alpha _{2i-2}}$}
\psfrag{m4}[tl]{$\scriptstyle{\alpha_{2i}}$}
\psfrag{m3p}[tl]{$\scriptstyle{\alpha' _{2i-1}}$}
\psfrag{n3}[tl]{$\scriptstyle{\beta _{i+1}}$}
\includegraphics[width=0.7\textwidth]{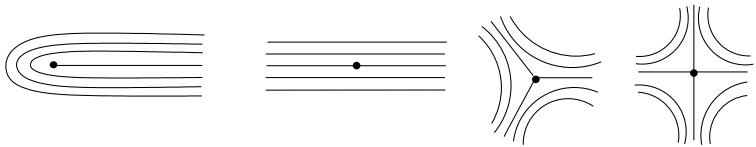}
\caption{Leaves near $p$--pronged singularities}\label{foliationtripsurvey}
\end{figure}

$\mu$ is called a \emph{transverse measure} on  $\mathcal{F}$. Also the singularities are classified with their number of prongs $p\geq 1$ ( a \emph{prong} is a piece of a leaf beginning at a singularity).  That is, at a $p$--pronged singularity the leaves of the foliation locally looks like as depicted in Figure \ref{foliationtripsurvey}.

A \textit{Whitehead move} on a measured foliation  contracts a compact leaf that joins two singularities or does the inverse.  Two measured foliations are \textit{Whitehead Equivalent} if one can be turned into the other after a finite number of Whitehead moves.   Let $\MF(M)$ denote the set of measured foliations on~$M$~up to isotopy and Whitehead equivalence. Let $(\mathcal{F},\mu)\in \MF(M)$ and  $k>0$.  Then $(\mathcal{F},k\mu)$ is a measured foliation with the same leaves as those of $(\mathcal{F},\mu)$ such that any arc $\alpha$ transverse to $\mathcal{F}$ has measure $k\mu(\alpha)$.  The space of \textit{Projective Measured Foliations},~$\mathcal{PMF}(M)$, on $M$~is the quotient space of $\mathcal{MF}(M)$ modulo~$\mathcal{(F,\mu)}\sim(\mathcal{F},k\mu)$, $k>0$. $\mathcal{PMF}(M)$ is homeomorphic to a sphere with dimension $6g+2s+2b-7$, and forms the boundary of the Teichm\"uller space $T(M)$ of $M$ which is an open ball of dimension $6g+2s+2b-6$  \cite{W88, FLP79}.  The closure $\overline{\mathcal{T(M)}}$ of $T(M)$ is a closed ball on which the mapping class group $\MCG(M)$ acts continuously.  By Brouwer Fixed Point Theorem each isotopy class has a fixed point on $T(M)$, and the analysis of this fixed point yields the Nielsen-Thurston classification theorem. This famous  theorem states that any homeomorphism of $M$  is isotopic to a  finite order or pseudo-Anosov or reducible homeomorphism ~\cite{W88, FLP79}. 

 A homeomorphism~$f \colon M\to M$~is \emph{pseudo\,-Anosov} if there exists a transverse pair of measured foliations $(\cF^s,\mu^s)$~and~$(\cF^u,\mu^u)$~and a number~$\lambda>1$ (the dilatation)~such that
\begin{align*}
f(\cF^s,\mu^s)&=(\cF^s,(1/\lambda)\mu^s)\\*
f(\cF^u,\mu^u)&=(\cF^u,\lambda \mu^u).
\end{align*}
We say that~$(\cF^s,\mu^s)$~and~$(\cF^u,\mu^u)$~are the \emph{stable} and \emph{unstable} invariant foliations of $f$ respectively.   A homeomorphism $f \colon M\to M$ is \emph{reducible} if it preserves a collection of mutually disjoint essential simple closed curves (\emph{reducing curves}). If some iterate of $f \colon M\to M$ is the identity, it is called \emph{finite order}. Each isotopy class is represented by a homeomorphism which is of one of these three types, and the isotopy class is named by the type of the homeomorphism it contains.  A pseudo-Anosov homeomorphism  in a pseudo-Anosov isotopy class is essentially unique. That is, if ~$g:M\to M$ is a pseudo\,-Anosov  homeomorphism in the isotopy class of $f$, then $f$ and $g$ are topologically conjugate by a homeomorphism that is isotopic to the identity. Let $[f]$ be a pseudo--Anosov isotopy class.  The induced action of  $[f]$ on $\overline{\mathcal{T(M)}}$  has exactly two fixed points which both lie on $\mathcal{PMF}(M)$, the projective classes $[\cF^s,\mu^s]$~and~$[\cF^u,\mu^u]$ of its  invariant measured foliations $(\cF^s,\mu^s)$~and~$(\cF^u,\mu^u)$,  and  every other point on $\mathcal{PMF}(M)$   converges to $[\cF^u,\mu^u]$ rapidly under the action of $[f]$ \cite{W88, FLP79}.  Since $\mathcal{MF}(M)$ has $\PIL$ (piecewise--integral--linear)  and  $\mathcal{PMF}(M)$ has $\PIP$ (piecewise--integral--projective) structure \cite{penner, W88, FLP79},  the action of $[f]$ on $\mathcal{PMF}(M)$ is piecewise projective linear and is locally described by integer matrices.

		 The matrix which acts on a piece that contains $[\mathcal{F}^u,\mu^u]$ in its closure has an eigenvalue $\lambda>1$ since $[\mathcal{F}^u,\mu^u]$ is a fixed point on $\mathcal{PMF}(M)$. Therefore,  if we compute the action of $[f]$ on $\mathcal{PMF}(M)$ and find a matrix with an eigenvalue $\lambda>1$ with associated eigenvector contained in the relevant piece,  the eigenvector corresponds to $[\mathcal{F}^u,\mu^u]$ and $\lambda$ gives the dilatation. In the case where the surface is the $n$-punctured disk $D_n$, we call such pieces  \emph{Dynnikov regions}, and the corresponding matrices  \emph{Dynnikov matrices}.  We coordinatize  the space of measured foliations $\mathcal{MF}_n$  making use of the Dynnikov's coordinates, and describe the action of $B_n$ on $\mathcal{MF}_n$ in terms of Dynnikov coordinates using the update rules.

%
%
%
%
%
%
%
%

\section{Dynnikov matrices of pseudo--Anosov braids}

In this section we reinterpret the aforementioned results from Thurston's theory of surface homeomorphisms  \cite{FLP79, W88} in terms of Dynnikov coordinates, and define Dynnikov matrices  which  describe the action of a given pseudo-Anosov braid  near its unstable  invariant measured foliation. We first define Dynnikov coordinates for measured foliations: Consider the set $\cA_n$ of Dynnikov arcs $\alpha_i $ ($1\le i\le 2n-4$) and $\beta_i$ ($1\leq i\leq n-1$) in~$D_n$.  Let $\gamma\in \cA_n$.   Given a measured foliation $(\mathcal{F},\mu)$ on $D_n$, we can use the measure $\mu$ on the foliation to define the \emph{ length} $\mu(\gamma)$ as,
$$\mu(\gamma)=\sup\sum_{i=1}^k{\mu(\gamma_i)}$$
where the supremum is taken over all finite collections $\gamma_1,\ldots,\gamma_k$ of mutually disjoint subarcs of $\gamma$ that are transverse to $\mathcal{F}$. Now, let $[\gamma_i]$ denote the isotopy class of 
$\gamma_i $  (under isotopies through~$\cA_n$).  Then, $[\gamma]$ has measure \[\mu([\gamma])=\inf_{\xi\in[\gamma]}\mu(\xi)\] which is well defined on~$\MF_n$.

Therefore the {\em Dynnikov coordinate function} $\rho:\mathcal{MF}_n\to\R^{2n-4}\setminus\{0\}$ is defined by

\begin{eqnarray*}
\rho(\mathcal{F},\mu) = (a;\,b)=(a_1,\ldots,a_{n-2},\,b_1,\ldots,b_{n-2}),
\end{eqnarray*}
where for $1\le i\le n-2$

\begin{eqnarray}
a_i=\frac{\mu([\alpha_{2i}])-\mu([\alpha_{2i-1}])}{2}\qquad
\text{and} \qquad b_i=\frac{\mu([\beta_i])-\mu([\beta_{i+1}])}{2} \label{dynnikov}
\end{eqnarray}

		Let $\cS_n=\R^{2n-4}\setminus\{0\}$ denote the space of Dynnikov coordinates of measured foliations on $D_n$. The Dynnikov coordinate function $\rho:\MF_n\to \cS_n$ is a homeomorphism when $\mathcal{MF}_n$ is endowed with its usual topology \cite{W88, FLP79}. The bijection is as given in Theorem \ref{lem:dynninvert} \cite{paper1, paper3}. That is,  let $(a,b)\in \mathbb{R}^{2n-4}\backslash\{0\}$. Then $(a,b)$ is the Dynnikov coordinate of exactly one element $(\mathcal{F},\mu)\in \mathcal{MF}_n$, which has
\begin{align}
\mu([\beta_i])&=2\max_{1\leq k\leq n-2}\left[|a_k|+\max(b_k,0)+\sum_{j=1}^{k-1}b_j\right]-2\sum_{j=1}^{i-1}b_j \\
\mu([\alpha_i])&= \left\{ \begin{array}{ll}
         (-1)^ia_{\lceil i/2 \rceil}+\frac{\mu([\beta_{\lceil i/2\rceil}])}{2}& \mbox{if $b_i \geq 0$};\\
        (-1)^ia_{\lceil i/2\rceil}+\frac{\mu([\beta_{1+\lceil i/2\rceil}])}{2}& \mbox{if $b_i \leq 0$}   
         \end{array} \right. 
\end{align}

Projectivizing Dynnikov coordinates yields an explicit bijection between $\cS_n/\mathbb{R}^{+}$ and $S^{2n-5}\cong\mathcal{PMF}_n$. Let $\mathcal{PS}_n$ denote the space of projective Dynnikov coordinates and $\rho(\mathcal{F},\mu)=(a;\,b)$. We shall write $[a;\,b]\in \mathcal{PS}_n$ to denote the Dynnikov coordinates of the projective class $[\mathcal{F},\mu]$ on $\mathcal{PMF}_n$. The Nielsen\,-Thurston Classification Theorem can be restated in terms of Dynnikov coordinates  \cite{paper1, paper3} as follows:

\begin{thm}

Let $\beta\in B_n$. Then

\begin{enumerate}[i.]
\item Either $\beta$ is of finite order, i.e. there is some $N>1$ such that $\beta^N(a;\,b)=~(a;\,b)$ for all~$(a;\,b)\in\cS_n$; or

\item $\beta$ is  reducible, i.e. there is some $(a;\,b)\in\mathcal{C}_n$ with $\beta(a;\,b)=(a;\,b)$; or
\item   $\beta$ is pseudo\,-Anosov, i.e. there is some~$(a^u;\,b^u)\in~\cS_n$ and a number $\lambda>1$ (the dilatation) such that
$\beta(a^u;\,b^u)=\lambda(a^u;\,b^u)$. In this case there is also some $(a^s;\,b^s)\in~\cS_n$ such that $\beta(a^s;\,b^s)=\frac{1}{\lambda}(a^s;\,b^s).$
\end{enumerate} 


\end{thm}

Let $\beta\in B_n$ be a pseudo--Anosov braid with unstable and stable invariant foliations having Dynnikov coordinates $(a^u;\,b^u)$ and $(a^s;\,b^s)$ respectively. Let $[a^u;\,b^u]$ and $[a^s;\,b^s]$ denote the projective classes of $(a^u;\,b^u)$ and $(a^s;\,b^s)$ on $\mathcal{PS}_n$ respectively.  Since the only fixed points of a pseudo--Anosov braid  are the projective classes of its invariant measured foliations  on $\mathcal{PMF}(M)$,  and  every other point on $\mathcal{PMF}(M)$   converges to the projective class of its  unstable invariant measured foliation under the action of $\beta$~\cite{W88, FLP79}, we can restate the following two lemmas  \cite{paper1, paper3} in terms of Dynnikov coordinates as follows:

\begin{lem}\label{fixedpointteichmullerdyn}
Let $\beta$ be a pseudo--Anosov braid. Then, the only fixed points of ~$\beta$~on
$\mathcal{PS}_n$ are $[a^u;\,b^u]$ and $[a^s;\,b^s]$.
\end{lem}
Therefore, any $(a;\,b)\in\cS_n$ satisfying $\beta(a;\,b)=k(a;\,b)$ for some $k>0$ is a multiple either of
$(a^u;\,b^u)$ or of $(a^s;\,b^s)$.

\begin{lem}
\label{lem:convergence2dynn}
Let $\beta$ be a pseudo--Anosov braid with fixed points $[a^u;\,b^u]$ and $[a^s;\,b^s]$ on $\mathcal{PS}_n$.  Then, for any~$[a;\,b]\in \mathcal{PS}_n$,~with~$[a;\,b]\neq [a^s;\,b^s]$,
\begin{align*}
\lim_{n \to \infty} \beta^n ([a;\,b])=[a^u;\,b^u].
\end{align*}
\end{lem}

\noindent Given an $n$--braid the update rules describe the action of $\beta$ on $\mathcal{PS}_n$. The collection of linear equations in various maxima in the update rules induce a piecewise linear action of $\beta$ on $\mathcal{PS}_n$, and hence  each region on $\mathcal{PS}_n$  has the structure of a polyhedron since it is described as a solution set for a system of linear inequalities induced by these equations. This is illustrated in Example \ref{1-2} where the piecewise linear action of the pseudo--Anosov $3$--braid $\sigma_1\sigma^{-1}_2$ on $\mathcal{PMF}_3\cong S^1$  is as given in Figure \ref{fig:teichmuller1}.

\begin{example}\label{1-2}

Let $(a;\,b)\in \mathcal{S}_3$, $\sigma_1(a;\,b)=(a';\,b')$ and $\sigma^{-1}_2(a';\,b')=(a'';\,b'')$. From Theorem~\ref{thm:update}

\begin{align*}
a'&=\tropical{\frac{ab}{a+1+b}} & b'&=\tropical{\frac{1+b}{a}}\\*
\intertext{that is}
a'&=a+b-\max\{a,0,b\} & b'&=\max\{0,b\}-a
\intertext{and}\\*
a''&=\tropical{\frac{a'}{a'b'+ 1+b'}} & b'&=\tropical{\frac{a'b'}{1+b'}}\\*
\intertext{that is}
a''&=a'-\max\{a'+b',0,b'\},& b''&=a'+b'-\max\{0,b'\}.
\end{align*}

Clearly, there are four main cases. Let us consider the  case $a\leq 0$,~$b\leq 0 $. We have

\begin{align*}
a'&=a+b ~~\text{and}~~ b'=-a
\end{align*}

\begin{align*}
a''&=a+b-\max(a+b-a,0,-a)& b''&=a+b-a-\max(0,-a)\\*
&=2a+b            &                         &=a+b
\end{align*}
Hence the action is given by the matrix;
\begin{align*}
D=\left[
\begin{array}{cc}
  2& 1  \\
   1& 1
\end{array}
\right].
\end{align*}
\begin{figure}[h!]
\begin{center}
\psfrag{2847}[tl]{$\scriptstyle{2.84}$} 
\psfrag{3847}[tl]{$\scriptstyle{3.84}$} 
\psfrag{1847}[tl]{$\scriptstyle{1.84}$} 
\psfrag{8009}[tl]{$\scriptstyle{8.00}$} 
\psfrag{2081}[tl]{$\scriptstyle{2.08}$} 
\psfrag{16669}[tl]{$\scriptstyle{16.67}$} 
\psfrag{4330}[tl]{$\scriptstyle{4.33}$} 
\psfrag{36540}[tl]{$\scriptstyle{34.70}$} 
\psfrag{16.90}[tl]{$\scriptstyle{16.90}$} 
\psfrag{70340}[tl]{$\scriptstyle{70.34}$} 
\psfrag{9012}[tl]{$\scriptstyle{9.01}$} 
\psfrag{R4}[tl]{$\scriptstyle{\mathcal{R}_4}$} 
\psfrag{R2}[tl]{$\scriptstyle{\mathcal{R}_2}$} 
\psfrag{R3}[tl]{$\scriptstyle{\mathcal{R}_3}$} 
\includegraphics[width=0.5\textwidth]{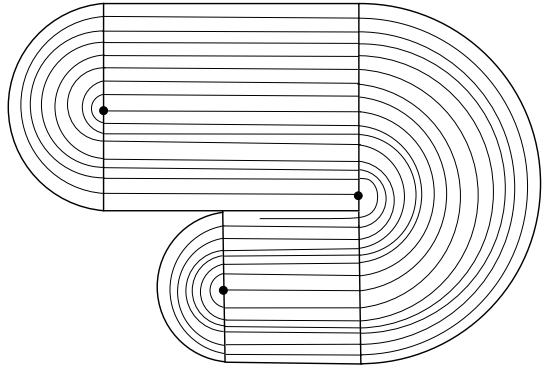}
\caption{The unstable invariant measured foliation for $\sigma_1\sigma^{-1}_2$} \label{1-2foli}
\end{center}
\end{figure}

\noindent which has an eigenvalue $\lambda=\frac{3+\sqrt{5}}{2}$ with corresponding eigenvector $p^u\approx -(0.850, 0.525)$. Since $p^u$ belongs to the region $a\leq 0$,~$b\leq 0 $, it is a fixed point for  $\sigma_1\sigma_2^{-1}$ on $\mathcal{PS}_3$. Hence, $p^u=[a^u;\,b^u]$ corresponds to the Dynnikov coordinates of the   unstable invariant measured foliation $[\mathcal{F}^u,\mu^u]$ by Lemma \ref{fixedpointteichmullerdyn}.  We can  compute  the measures assigned by $(\mathcal{F}^u,\mu^u)$ on the Dynnikov arcs  from its Dynnikov coordinates $(-0.850, -0.525)$ by $(6)$ and $(7)$.  Therefore, the  approximate measures on the Dynnikov arcs are given by $\mu(\alpha_{1})=2.226,~  \mu(\alpha_{2})=0.525,~  \mu(\beta_{1})= 1.70,  ~\mu(\beta_{2})=2.751$.  Replacing each Dynnikov arc  with a Euclidean rectangle of height given by these measures  and length 1,  and  endowing each rectangle with a horizontal measured foliation (the transverse measure is induced from the Euclidean metrics on the rectangles),  and gluing the vertical sides of the rectangles we can reconstruct $(\mathcal{F}^u,\mu^u)$ depicted in  Figure \ref{1-2foli}. The construction is analogous to the one in Figure \ref{const1}. We note that  the horizontal leaves in the rectangles are glued together in a unique measure preserving way since the triangle inequality is satisfied in each triangular region.

Similarly, if $a\geq 0,~b\geq 0,~ b\leq a\leq 2b$ we compute  the matrix
$$
\left[ \begin {array}{cc} 1&-1\\\noalign{\medskip}-1&2
\end {array} \right]
$$
which has an eigenvalue $1/\lambda$ such that the associated eigenvector $p^s$ belongs to the region $a\geq 0,~b\geq 0,~ b\leq a\leq 2b$. Hence $p^s$ is a fixed point and corresponds to the   stable invariant measured foliation $[\mathcal{F}^s,\mu^s]$. Other matrices are computed similarly and the action of  $\sigma_1\sigma^{-1}_2$ on $\mathcal{PMF}_3$ in terms of Dynnikov coordinates is illustrated in Figure \ref{fig:teichmuller1}.

\begin{figure}[h!]
\begin{center}
\psfrag{2}[tl]{$\tiny{\left[ \begin {array}{cc} 1&-1\\\noalign{\medskip}1&0
\end {array} \right]}$}
\psfrag{u}[tl]{$\tiny{p^u}$}
\psfrag{s}[tl]{$\tiny{p^s}$}
\psfrag{A}[tl]{$\tiny{A}$}
\psfrag{B}[tl]{$\tiny{B}$}
\psfrag{a}[tl]{$\tiny{a}$}
\psfrag{b}[tl]{$\tiny{b}$}
\psfrag{C}[tl]{$\tiny{C}$}
\psfrag{D}[tl]{$\tiny{D}$}
\psfrag{E}[tl]{$\tiny{E}$}
\psfrag{l1}[tl]{$\tiny{\ell_1}$}
\psfrag{l2}[tl]{$\tiny{\ell_2}$}
\psfrag{F}[tl]{$\tiny{F}$}
\psfrag{3}[tl]{$\tiny{\left[ \begin {array}{cc} 1&-1\\\noalign{\medskip}-1&2
\end {array} \right]}$}
\psfrag{4}[tl]{$\tiny{\left[ \begin {array}{cc} 0&1\\\noalign{\medskip}-1&2
\end {array} \right]}$}
\psfrag{1}[tl]{$\tiny{\left[ \begin {array}{cc} 2&-1\\\noalign{\medskip}1&0
\end {array} \right]}$}
\psfrag{5}[tl]{$\tiny{\left[ \begin {array}{cc} 0&1\\\noalign{\medskip}-1&1
\end {array} \right]}$}
\psfrag{6}[tl]{$\tiny{\left[ \begin {array}{cc} 2&1\\\noalign{\medskip}1&1
\end {array} \right]}$}
\psfrag{a=b}[tl]{${a=b}$}
\psfrag{a=2b}[tl]{${a=2b}$}
\includegraphics[width=0.9\textwidth]{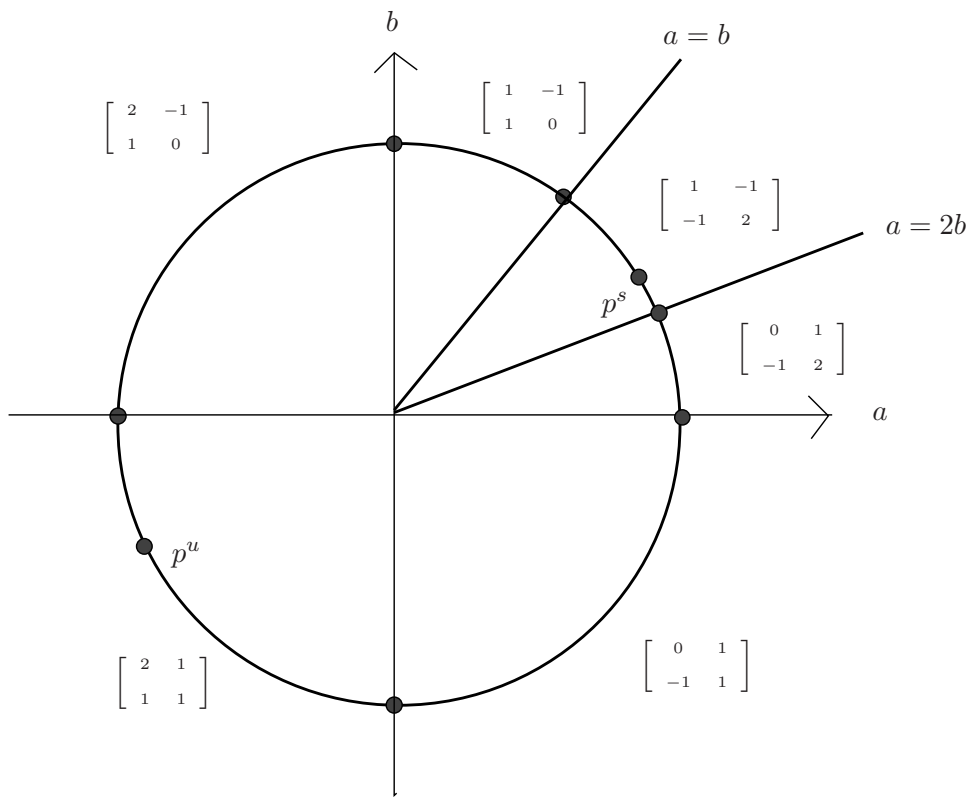}
\caption{The action of $\sigma_1\sigma^{-1}_2$ on $\mathcal{PS}_3$}
\label{fig:teichmuller1}
\end{center}
\end{figure}

\end{example}
In this example, $D=\left[
\begin{array}{cc}
  2& 1  \\
   1& 1
\end{array}
\right]$ is the  Dynnikov matrix for $\sigma_1\sigma^{-1}_2$. Given a pseudo--Anosov $n$--braid $\beta$ a Dynnikov matrix for $\beta$ is defined as follows: Let $(\mathcal{F},\mu)$ be the  unstable invariant  measured foliation of $\beta$ with Dynnikov coordinates $(a^u;\,b^u)$. Let $[a^u;\,b^u]$ denote the projective class of $(a^u;\,b^u)$ on $\mathcal{PS}_n$. Each piece (with respect to the piecewise projective action) $\mathcal{R}_i\subset \mathcal{PS}_n$ that contains $[a^u;\,b^u]$ in its closure is called a \emph{Dynnikov region}. A \emph{Dynnikov matrix}  $D_i:\mathcal{R}_i\to\mathcal{PS}_n$, ($1\leq i \leq k$)  for $\beta$ is a $(2n-4)\times(2n-4)$ integer matrix which describes the action of $\beta$ on $\mathcal{R}_i$.

\noindent We note that there can be more than one Dynnikov region  for a given pseudo--Anosov braid if $p^u=[a^u, b^u]$  is on the boundary of several regions on $\mathcal{PS}_n$.   Consider  for example the $6$--braid $\beta=\sigma_1\sigma_2\sigma_3\sigma_4\sigma^{-1}_5$. Resolving the various maxima in the update rules as in Example \ref{1-2} for the region  $a_i\leq 0, b_i\leq 0 $ ($1\leq i \leq 4$) gives four matrices. Each of these matrices has an eigenvalue $\lambda\approx 2.081$ such that the  associated eigenvector  $(a; b)\in \mathbb{R}^8\setminus\left\{0\right\}$ is given by

\begin{align*}
a_1&\simeq -1& b_1&\simeq -2.081\\*
a_2&\simeq -3.081&b_2&\simeq -4.330\\*
a_3&\simeq -7.411&b_3&\simeq -9.012\\*
a_4&\simeq -18.27&b_4&\simeq -16.904.
\end{align*}

 \noindent $p^u$ satisfies the inequalities $a_i\leq 0, b_i\leq  0$ ($1\leq i \leq 4$) and the equalities $a_2-a_1=b_1$, $a_3-a_2=b_2$.  Therefore, each of these matrices is a Dynnikov matrix for $\beta$.  An important observation regarding the singularity structure of $(\mathcal{F},\mu)$ is as follows: 
By Theorem \ref{lem:dynninvert} one can work out $\alpha_j$ for all $1\leq j\leq 8$ and $\beta_j$ for all 
$1\leq j\leq 5$, and observe that $|\alpha_2|-b_1=|\alpha_4|-b_2=|\alpha_6|-b_3$. This implies that there exists a leaf which joins three $3$-pronged singularities, and a Whitehead move contracts this leaf and yields a $5$ pronged singularity. This leads to the following natural question:

\begin{problem1}
Let $D$ be a Dynnikov matrix for  a pseudo--Anosov braid $\beta\in B_n$ with unstable invariant  measured foliation $(\mathcal{F}^u,\mu^u)$. Give an algorithmic way to determine the singularity structure of $(\mathcal{F}^u,\mu^u)$
  from $D$.  \end{problem1}

 In Example \ref{1-2} we computed  the action of $\beta=\sigma_1\sigma^{-1}_2$ on the whole space $\mathcal{PMF}_3$. However, this is not necessary to compute Dynnikov matrices. To be more specific,  since  $[\mathcal{F}^u,\mu^u]$ is a globally attracting fixed point for the action of $\beta$  on $\mathcal{PMF}_n$ (Lemma \ref{lem:convergence2dynn}) it is easy to find Dynnikov regions and hence Dynnikov matrices as will be explained in more detail in Section \ref{topent}. Indeed,  observe  in Example \ref{1-2} that each region on  $\mathcal{PMF}_3$ is attracted to the third quadrant  which is the Dynnikov region for $\beta$. If $\ell_1$ and $\ell_2$ denote the lines $a=b$ and $a=2b$ in the first quadrant, and $+a$~($-a$) and $+b~(-b)$ denote the positive (negative) $a$-axis and $b$-axis respectively, we get  $\ell_1\rightarrow+b\rightarrow-a$, $\ell_2\rightarrow+a\rightarrow-b$ where $x\rightarrow y$ means $x$ is sent onto $y$ by the action of $\beta$.

\section{Topological entropy of infinite families of pseudo-Anosov braids}\label{topent}
In \cite{paper1} Dynnikov matrices were used to introduce a fast method for computing the topological entropy of each member of an infinite family of pseudo--Anosov braids. The topological entropy of an isotopy class is defined to be the minimum topological entropy of a homeomorphism contained in it. When the isotopy class is pseudo--Anosov with dilatation $\lambda$, all pseudo--Anosovs in the class have the minimum topological entropy which is $\log \lambda $. Thus, the topological entropy of a pseudo--Anosov isotopy class equals  $\log \lambda $ \cite{FLP79}.  The usual approach to compute the topological entropy of an isotopy class of surface homeomorphims is to use train-track methods such as the Bestvina--Handel algorithm \cite{BH95}. If the isotopy class is pseudo--Anosov, the algorithm gives  a train track associated with a transition matrix from which a Markov partition for the pseudo-Anosov homeomorphism in the isotopy class is constructed, and hence the dilatation, the invariant foliations and the minimal periodic orbit structure of the isotopy class are obtained. 

	In \cite{M06} Moussafir introduced an alternative method for computing the topological entropy of braids making use of Dynnikov coordinates and the update rules.  The main idea of his approach lies in the following result of Thurston \cite{W88, FLP79}:  If $\lambda$ is the dilatation of a given pseudo\,-Anosov  braid $\gamma\in B_n$, then for any essential simple closed curves $\alpha$ and $\beta$, the geometric intersection number $i(\gamma^n(\alpha),\beta)$ grows like 
$C\times\lambda^n$ as $n\to \infty$ for some positive constant $C\in \mathbb{R}$. More precisely, he constructs for every braid $\gamma\in B_n$, an integer sequence
$$c_m=\frac{1}{m}\log c(\gamma^m \cL)$$ where  $\cL$ is a particular type of integral lamination with Dynnikov coordinates $\rho(\cL)=(a;\,b)$, $\rho(\gamma^m \cL)$ is obtained from the update rules and $c(\gamma^m \cL)$ denotes the minimum number of intersections of $\gamma^m \cL$ with the $x$-axis. The growth rate of $c_m$ gives an estimate for the topological entropy of $\beta$. The major advantage of his method is that it works much faster and is more direct than the classical train-track approach \cite{BH95}. However, the method is numerical and only gives an estimate for braid entropy. On the other hand, the fact that each Dynnikov matrix has an eigenvalue $\lambda>1$ associated with the Dynnikov coordinates $[a^u;\,b^u]$ of $[\mathcal{F}^u,\mu^u]$, and that $\lambda$ gives the dilatation of $\gamma$  yields the method introduced in \cite{paper1, paper3}   to compute the exact topological entropy. 
Next, we discuss how we use  Dynnikov matrices to compute the topological entropy of each member of a given pseudo--Anosov braid family. We shall consider the  two families of braids studied in \cite{paper1}. These families are 

\begin{align*}
\beta_{m,n} &= \sigma_1\sigma_2\cdots \sigma_{m-1}\sigma_m~ \sigma_{m+1}\I\sigma_{m+2}\I\cdots
\sigma_{m+n-1}\I\sigma_{m+n}\I \in
B_{m+n+1}, \quad\text{and}\\*
\sigma_{m,n} &= \sigma_1\sigma_2\cdots \sigma_{m-1}\sigma_m~ \sigma_m\sigma_{m-1}\cdots\sigma_2\sigma_1~
\sigma_1\sigma_2\cdots \sigma_{m+n-1}\sigma_{m+n} \in B_{m+n+1}.
\end{align*}

To compute the topological entropy of each member of a family of braids using Dynnikov matrices \cite{paper1},  we compute Dynnikov regions and matrices for enough braids in the family  until a general pattern is spotted, conjecture that the pattern holds for all braids in the family, and then prove the conjecture \cite{paper1}. Next, we shall briefly explain the  method  on the subfamily $\tau_n=\beta_{n-2,1}=\sigma_1\sigma_2\cdots\sigma_{n-2}\sigma^{-1}_{n-1}\in B_n$. That is, the following steps give a recipe to find the Dynnikov coordinates of the invariant foliation and a Dynnikov matrix for each braid $\tau_n$. 
\begin{itemize}

 \item \textbf{Step 1:} (Experiment) Since $[\mathcal{F}^u,\mu^u]$ is a globally attracting fixed point for the action of $\tau_{n}\in B_n$ on $\mathcal{PMF}_n$, it is easy to find its Dynnikov coordinates numerically. We use the program {\tt Dynn} \cite{Toby}  for this. The program picks a random point $(a;\,b)\in \mathbb{R}^{2n-4}\setminus\{0\}$ and iterates it with the given braid $\tau_n$ until it arrives in a region $\mathcal{R}$ in which there exists a point $(a^u;\,b^u)$ with $D[a^u;\,b^u]=[a^u;\,b^u]$, where $D$ describes the action of $\tau_n$ in $\mathcal{R}$.  Thus, $(a^u;\,b^u)$ corresponds to the Dynnikov coordinates of $[\mathcal{F},\mu]$ and $D$ is a Dynnikov matrix.  We obtain Dynnikov matrices of $\tau_n$ for different values of $n$. We observe that the Dynnikov coordinates of $(a^u;\,b^u)$ are all negative and hence decide to compute the update rules under the assumption that $a_j\leq 0$ and $b_j\leq 0$ for all $1\leq j \leq n-2$ as the second step.

 \item \textbf{Step 2:} We note that there can be more than one  Dynnikov matrix if $[\mathcal{F},\mu]$ is on the boundary of several Dynnikov regions. For example, it follows for $\tau_n$ that there are $2^{n-4}$ Dynnikov regions adjacent to $(a^u;\,b^u)$ since for each $1\leq j\leq n-4$ the update rules can be calculated either under the asumption that $a_{j+1}-a_j\leq b_j$ or under the assumption that $a_{j+1}-a_j\geq b_j$. We choose a region where $a_{j+1}-a_j\leq b_j$ for ~$1\leq j \leq n-3$ and write $\mathcal{R}^{(n)}$ to denote this region, compute the update rules for $\tau_n$ in $\mathcal{R}^{(n)}$ and write $D^{(n)}$ for the update matrix in $\mathcal{R}^{(n)}$.
 
  \item \textbf{Step 3:} We prove that for each $n$, $D^{(n)}$ is a Dynnikov matrix. To do this we find a general form for the characteristic polynomial $f_n(x)$ of $D^{(n)}$ for each $n$ and prove that $f_n(x)$ has eigenvalue $r>1$ with corresponding eigenvector $(a^u;\,b^u)$ contained in $\mathcal{R}^{(n)}$, and hence conclude that $\mathcal{R}^{(n)}$ is a Dynnikov region and $D^{(n)}$ is a Dynnikov matrix. For $\tau_n$,  $f_n(x)$ is given by
  
  \begin{eqnarray*}
f_n(x)=(x+1)^{n-4}(x^n-2x^{n-1}-2x+1)
\end{eqnarray*}

which has a root $r>1$, and the unique eigenvector corresponding to $r>1$ in $\mathcal{R}^{(n)}$ is given by
 
\begin{align*}
a_j &=
\begin{cases}
-r(r^j-1);&1\leq j <n-2 \\*
-(r^{n-1}-1)(r-1) & j=n-2,
\end{cases}
\\
b_j &=
\begin{cases}
-r^{j+1}(r-1) & 1\leq j <n-2\\*
-(r^{n-1}-1) & j=n-2.
\end{cases}
\end{align*}

\end{itemize}

Applying the same approach to $\beta_{m,n}$ and $\sigma_{m,n}$ yields the following theorems \cite{paper1}.

\begin{thm}
\label{thm:betamn}
Let $m,n\ge 1$. Then $\beta_{m,n}\in B_{m+n+1}$ is a pseudo\,-Anosov 
braid, whose dilatation~$r$ is the unique root in $(1,\infty)$ of the
polynomial
\[f_{m,n}(r)=(r-1)(r^{m+n+1}-1) - 2r(r^m+r^n).\]
The Dynnikov coordinates $(a^u;\,b^u)\in\cS_{m+n+1}$ of the unstable
invariant measured foliation of~$\beta_{m,n}$ are given by
\begin{align*}
a_i &= 
\begin{cases}
-r(r^n+1)(r^i-1) & \text{ if \ } 1\le
  i\le m-1 \\
-(r^{m+1}-1)(r^{n+1}-1) & \text{ if \  } i =
  m\\
-(r^{m+1}-1)(r^{m+n+1-i}-1)r^{i-m}
  & \text{ if \ } m+1 \le i \le m+n-1,
\end{cases}
\\
b_i &= 
\begin{cases}
-(r-1)(r^n+1)r^{i+1} & \text{ if \ } 1\le
  i\le m-1 \\
-(r+1)(r^{m+1}-1) & \text{ if \ } i =
  m\\
-(r-1)(r^{m+1}-1)r^{i-m} \text{\phantom{Some st\,}}
  & \text{ if \ } m+1 \le i \le m+n-1.
\end{cases}
\end{align*}
\end{thm}

\begin{thm}[The braids $\sigma_{m,n}$ for $n\ge m+2$]
\label{thm:sigmamn}
Let $1\le m\le n-2$. Then $\sigma_{m,n}\in B_{m+n+1}$ is a pseudo\,-Anosov 
braid, whose dilatation~$r$ is the unique root in $(1,\infty)$ of the
polynomial
\[g_{m,n}(r)=(r-1)(r^{m+n+1}+1) + 2r(r^m-r^n).\]
The Dynnikov coordinates~$(a^u;\,b^u)\in\cS_{m+n+1}$~of the unstable
invariant measured foliation of~$\sigma_{m,n}$~are given by
\begin{align*}
a_i &= 
\begin{cases}
(r^n-1)(r^{i+1}-1)r & \text{ if \ } 1\le i\le m-1\\*
(r^{m+1}-1)(r^{m+n-i}-1)r^{i+1-m} & \text{ if \ } m\le i \le m+n-1,
\end{cases}
\\
b_i &= 
\begin{cases}
(r-1)(r^n-1)r^{i+1} & \text{ if \ }1\le i\le m-1\\*
(r-1)(r^{m+1}-1)r^{i-m} \text{\phantom{Some st\,\,}}& \text{ if \ }m\le i\le m+n-1.
\end{cases}
\end{align*}
\end{thm}

\begin{remark}

 Dynnikov and train track transition matrices both yield a way to compute the topological entropy of pseudo--Anosov braids. However,  finding Dynnikov matrices is much easier than finding train track transition matrices.  See \cite{paper1, paper3} for  an experimental comparison of these matrices.  Indeed, a recent work on the Nielsen--Thurston classification  problem  \cite{MSY} has shown that such matrices can be computed in polynomial time. 
\end{remark}

\section{Dynnikov  and train track transition matrices}\label{topent2}
In \cite{b10} Birman, Brinkman and Kawamuro investigate the spectrum of the train track transition matrix $T$ of a given pseudo-Anosov homeomorphism $f$ on an orientable surface and show that the characteristic polynomial of $T$ factors into three polynomials: the first has the dilatation $\lambda$ of $f$ as its largest root; the second relates to the action of $f$ on the singularities of the invariant foliations $(\mathcal{F}^u,\mu^u)$ and $(\mathcal{F}^s,\mu^s)$; and the third relates to the degeneracies of a symplectic form introduced in \cite{penner}. The aim in \cite{paper3} is to show that any Dynnikov matrix shares the same set of eigenvalues with any train track transition matrix up to roots of unity and zeros --- this is proved in some cases but not in full generality. The properties of the spectrum of $T$, which is difficult to calculate, can therefore be studied using the spectrum of $D$, which is easy to calculate. The following is a brief summary of the results stated in \cite{paper3}.

	A \emph{train track} $\tau$ on $D_n$ is a one dimensional CW complex smoothly embedded on $D_n$ such that at each switch there is a unique tangent vector. The vertices are called \emph{switches}, and the edges are called \emph{branches}. We require that every component of $D_n-\tau$ is either a once-punctured $p$-gon with $p\geq 1$ or an unpunctured $k$-gon with $k\geq 3$ (where the boundary of $D_n$ is regarded as a puncture). A \emph{transverse measure} on $\tau$ assigns a measure to each branch of~$\tau$ such that at each switch of $\tau$ the measures satisfy the \emph{switch conditions} (some particular linear equations). 
	
	 \emph{Measured train tracks} are train tracks endowed with a transverse measure which provide coordinate patches for measured foliations and integral laminations \cite{mosher1, penner}.  We write $\mathcal{W}(\tau)$ for the space of transverse measures associated to $\tau$ and say that a foliation $(\mathcal{F},\mu)$ is \emph{carried} by $\tau$ if it arises from some non--negative transverse measure in $\mathcal{W}(\tau)$. In particular, there is a homeomorphism from the space of non--negative transverse measures $\mathcal{W}^{+}(\tau)$ on~$\tau$ to the space of measured foliations $\mathcal{MF}(\tau)$ carried by $\tau$ \cite{mosher1, penner, paper3}. 
 
 	Every pseudo\,-Anosov  homeomorphism $f$ has an \emph{invariant train track}: that is, a train track $\tau$ whose image under $[f]$ is another train track which can be collapsed onto $\tau$ in a regular neighbourhood of $\tau$ in a smooth way. The associated \emph{train track transition matrix} $T$ has entries $T_{i,j}$ given by the number of occurences of $e_i$ in the edge path $f(e_j)$ where $e_1,\dots,e_k$ denote the branches of $\tau$. The action of a given pseudo\,-Anosov  isotopy class $[f]$ on $\mathcal{W}(\tau)$ is given by the transition matrix associated to the invariant train track of $[f]$.   The largest eigenvalue of $T$ equals the dilatation $\lambda$ and the entries of the unique (up to scale) associated eigenvector $v^u$ in  $\mathcal{W}(\tau)$ are strictly positive and correspond to $(\mathcal{F}^u,\mu^u)$. 
	
	In \cite{paper3} it is shown that the isospectrality of Dynnikov  and train track transition matrices of a pseudo--Anosov braid depends on the singularity structure of the invariant foliations, and how the prongs of the singularities  are permuted under the action of the braid.  The proofs of following theorems  can be found in  \cite{paper3}. Each proof is based on the approach explained briefly under the relevant theorem.
	
\begin{thm}\label{firstthm}
Let $\beta\in B_n$ be a pseudo\,-Anosov braid with unstable invariant measured foliation $(\mathcal{F}^u,\mu^u)$. Let  $\tau$ be an invariant train track with associated transition matrix $T$. If $(\mathcal{F}^u,\mu^u)$ has only unpunctured $3$-pronged and punctured $1$-pronged singularities, then $\beta$ has a unique Dynnikov matrix $D$, and $D$ and $T$ are isospectral.
\end{thm}

If $(\mathcal{F}^u,\mu^u)$  has only unpunctured $3$-pronged and punctured $1$-pronged singularities (which is a generic property in $\mathcal{PMF}_n$), it is carried by a \emph{complete} train track $\tau$. That is, $\tau$ has the property that each component of $D_n -\tau$ is either a trigon or a once punctured monogon. In this case, $\mathcal{MF}(\tau)$ defines a chart in $\mathcal{MF}_n$ containing $(\mathcal{F}^u,\mu^u)$ in its interior and hence there is a unique Dynnikov matrix $D$. We use the change of coordinate function $\mathcal{W}^+(\tau)\to\cS_n$ which is linear in some neighbourhood of  $v^u$  (the unstable invariant measured foliation in train track coordinates) in $\mathcal{W}^+(\tau)$,  and show that $D$ and $T$ share the same spectrum.

\begin{thm}\label{thm1}
Let $\beta\in B_n$ be a pseudo\,-Anosov braid with unstable invariant measured foliation $(\mathcal{F}^u,\mu^u)$ and dilatation $\lambda>1$. Let $\tau$  be an invariant train track of $\beta$ with associated transition matrix $T$. If $\beta$ fixes the prongs at all singularities other than unpunctured $3$-pronged and punctured $1$-pronged singularities, then any Dynnikov matrix $D_i$ is isospectral to $T$  up to some eigenvalues~$1$.
\end{thm}

\begin{thm}\label{uniquematrix}
Let $\beta\in B_n$ be a pseudo\,-Anosov braid with unstable invariant measured foliation $(\mathcal{F}^u,\mu^u)$ and dilatation $\lambda>1$. Let $\tau$  be an invariant train track of $\beta$ with associated transition matrix $T$. If all components of $D_n-\tau$ are odd-gons and $\beta$ fixes the prongs at all singularities other than unpunctured $3$-pronged and punctured $1$-pronged singularities, then there is a unique Dynnikov region.
\end{thm}

 If $(\mathcal{F}^u,\mu^u)$ has singularities other than unpunctured $3$-pronged and punctured $1$-pronged singularities, then $\tau$ is not complete and therefore does not define a chart. There are two subcases to consider: first, when $\beta$ fixes the prongs of $(\mathcal{F}^u,\mu^u)$; and second, when it permutes them non-trivially. Whichever is the case, we construct a complete train track $\tau_p$ from $\tau$ so that $(\mathcal{F}^u,\mu^u)$ is carried by $\tau_p$ and is contained in the interior of $\mathcal{\MF}(\tau_p)$. In order to do this, we make use of the \emph{pinching move} \cite{mosher1, penner}. However, $\tau_p$ is not an invariant train track  unless relevant prongs of $(\mathcal{F}^u,\mu^u)$ are fixed by $\beta$. Therefore, we use another move called the \emph{diagonal extension move} \cite{mosher1, penner} which constructs several complete train tracks that give a set of charts which fit nicely in $\mathcal{MF}(\tau_p)$, with the property that the action in each of them is described explicitly. The results stated above  are mainly based on the interplay between the charts constructed from these two different moves.  The claim stated in Problem 2 has been observed in a wide range of examples \cite{Toby} but has not yet been proven.

\begin{problem2}
Let $\beta\in B_n$ be a pseudo\,-Anosov  braid  with train track transition matrix $T$ and Dynnikov matrices $D_i$. If $\beta$ permutes the prongs of its invariant unstable measured foliation non-trivially, is every Dynnikov matrix $D_i$ isospectral to $T$ up to roots of unity? 
\end{problem2}

We encourage the reader to study for example the $6$-braid $\beta=\sigma_1\sigma_2\sigma_3\sigma_4\sigma^{-1}_5$ on $D_6$ with invariant unstable measured foliation $(\mathcal{F}^u,\mu^u)$ having six 1--pronged, and one $5$--pronged singularities whose prongs are permuted non-trivially under the action of the braid. $\beta$ has four Dynnikov matrices acting in Dynnikov regions defined by the inequalities $a_i\leq 0, b_i\leq  0$ ($1\leq i \leq 4$) and the equalities $a_2-a_1=b_1$, $a_3-a_2=b_2$, and each Dynnikov matrix is isospectral to $T$ up to roots of unity. 
 
%


\end{document}